\documentclass[12pt]{article}
\usepackage{amsmath,amssymb,amsthm,amscd}
\usepackage{color}
\usepackage[T1]{fontenc}
\RequirePackage{easymat}
\DeclareMathOperator{\rank}{rank}

\def\proof{{\noindent \it Proof. }}
\theoremstyle{definition}

\theoremstyle{remark}

\theoremstyle{plain}
\newtheorem{theorem}{Theorem}

\newtheorem{lemma}[theorem]{Lemma}

\newtheorem{example}[theorem]{Example}

\newtheorem {algorithm} [theorem] {Algorithm}

\renewcommand{\ge}{\geqslant}
\renewcommand{\le}{\leqslant}

\title{Regularizing algorithm for
mixed matrix pencils}
\author{Tetiana Klymchuk\\ Departament de Matem\`{a}tiques,\\ Universitat Polit\'{e}cnica de Catalunya\\
Barcelona, SPAIN\\tetiana.klymchuk@upc.edu}

\begin{document}
 \maketitle

\abstract{
\normalsize \vspace*{2mm} {\small
P. Van Dooren (1979) constructed an algorithm for computing all singular summands of Kronecker's canonical form of
a matrix pencil. His algorithm uses only unitary transformations, which improves its numerical stability. We extend Van
Dooren's algorithm to square complex matrices with respect to consimilarity transformations $A\mapsto
SA\bar S^{-1}$ and to pairs of $m\times
n$ matrices up to transformations
$(A,B)\mapsto (SAR,SB\bar{R})$, in
which $S$ and $R$ are nonsingular
matrices.
}
}

\newcommand{\p}{\underline}
\newcommand{\ca}{\mathcal}
\renewcommand{\baselinestretch}{1}
\normalsize

\section{Introduction}

Van Dooren \cite{van} gave an algorithm
that for each pair $(A,B)$ of complex
matrices of the same size constructs
its \emph{regularizing decomposition};
that is, it constructs a matrix pair
that is simultaneously equivalent to
$(A,B)$ and has the form
\[
(A_1,B_1)\oplus\dots\oplus (A_t,B_t)
\oplus (\p A,\p B)
 \]
in which $(\p A,\p B)$ is a pair of
nonsingular matrices and each other
summand has one of the forms:
\[
(F_n,G_n),\quad (F_n^T,G_n^T),\quad (I_n,J_n(0)),\quad (J_n(0),I_n),
\]
where $J_n(0)$ is the singular Jordan
block and
\[
F_n:=\begin{bmatrix}
0&&0\\[-6pt]1&\ddots&\\[-6pt]&\ddots&0\\0&&1
\end{bmatrix},\quad
G_n:=\begin{bmatrix}
1&&0\\[-6pt]0&\ddots&\\[-6pt]&\ddots&1\\0&&0
\end{bmatrix}
\]
are $n\times(n-1)$ matrices; $n\ge 1$.
Note that $(F_1,G_1)=(0_{10},0_{10})$;
we denote by $0_{mn}$ the zero matrix
of size $m\times n$, where
$m,n\in\{0,1,2,\dots\}$. The algorithm
uses only unitary transformations,
which improves its computational
stability.

We extend Van Dooren's algorithm to
square complex matrices up to
consimilarity transformations $A\mapsto
SA\bar S^{-1}$ and to pairs of $m\times
n$ matrices up to transformations
$(A,B)\mapsto (SAR,SB\bar{R})$, in
which $S$ and $R$ are nonsingular
matrices.

A regularizing algorithm for matrices
of undirected cycles of linear mappings
was constructed by Sergeichuk
\cite{ser} and, independently, by Varga
\cite{var}. A regularizing algorithm
for matrices under congruence was
constructed by Horn and Sergeichuk
\cite{h-s}.

All matrices that we consider are
complex matrices.

\subsection{Regularizing unitary algorithm for matrices under consimilarity}

\quad\\

Two matrices $A$ and $B$ are
\emph{consimilar} if there exists a
nonsingular matrix $S$ such that $
SA\bar S^{-1}=B. $ Two matrices are
consimilar if and only if they give the
same semilinear operator, but in
different bases. Recall that a mapping
$\mathcal A:U\to V$ between complex
vector spaces is \emph{semilinear} if
\[
\mathcal A(au_1+bu_2)=\bar a\mathcal Au_1+\bar b\mathcal Au_2
\]
for all $a,b\in\mathbb C$ and
$u_1,u_2\in U$.

The canonical form of a matrix under
consimilarity is the following (see
\cite{hon-hor} or \cite{j-h}): \\
{\it
Each square complex matrix is
consimilar to a direct sum, uniquely
determined up to permutation of direct
summands, of matrices of the following
types:
\begin{itemize}
  \item a Jordan block
      ${J_k(\lambda)}$ with
      ${\lambda\ge 0}$, and
  \item $\begin{bmatrix} 0&1\\\mu
      &0\end{bmatrix}$ with
      ${\mu\notin\mathbb R}$ or
      ${\mu<0}$.
\end{itemize}}

Thus, each square matrix $A$ is
consimilar to a direct sum
\[
J_{n_1}(0)\oplus\dots\oplus J_{n_k}(0)\oplus \underline A,
\]
in which $\underline A$ is nonsingular
and is determined up to consimilarity;
the other summands are uniquely
determined up to permutation. This sum
is called a \emph{regularizing
decomposition} of $A$. The following
algorithm admits to construct a
regularizing decomposition using only
unitary transformations.

\begin{algorithm}
{\it Let $A$ be a singular $n\times n$
matrix. By unitary transformations of
rows, we reduce it to the form
\[ S_1A=
\begin{bmatrix}
  0_{r_1n} \\
  A' \\
\end{bmatrix},
\quad S_1 \text{ is unitary, }
\]
in which the rows of $A'$ are linearly
independent. Then we make the
coninverse transformations of columns
and obtain
\begin{equation*}\label{A1}
S_1A\bar {S_1}^{-1}=
\begin{bmatrix}
  0_{r_1} & 0 \\
  \star & A_1 \\
\end{bmatrix}
\end{equation*}
We apply the same procedure to $A_1$
and obtain
\begin{equation*}\label{A2}
S_2A_1\bar{S_2}^{-1}=
\begin{bmatrix}
  0_{r_2} & 0 \\
  \star & A_2 \\
\end{bmatrix},
\quad S_2 \text{ is unitary, }
\end{equation*}
in which the rows of $[\star \ A_2]$
are linearly independent.

We repeat this procedure until we
obtain
\begin{equation*}\label{At}
{S}_tA_{t-1}\bar{S_t}^{-1}=
\begin{bmatrix}
  0_{r_t} & 0 \\
  \star & A_t \\
\end{bmatrix},
\quad S_t \text{ is unitary, }
\end{equation*}
in which $A_t$ is nonsingular.  The
result of the algorithm is the sequence
$ r_1,r_2,\dots,r_{t},A_t.$ }
\end{algorithm}

For a matrix $A$ and a nonnegative
integer $n$, we write
\[
A^{(n)}:=
           \begin{cases}
             0_{00}, & \hbox{if $n=0$;} \\
A\oplus\dots\oplus A\ (n\text{ summands}),
 & \hbox{if $n\ge 1$.}
           \end{cases}
\]
\begin{theorem}\label{thrm1}
 Let $r_1, r_2, \dots, r_{t}, A_t$ be obtained by applying Algorithm 1 to a square complex matrix $A$. Then
 \begin{equation}\label{jjh}
 r_1\ge r_2\ge \dots\ge r_{t}
 \end{equation}
 and
 $A$ is consimilar to
 \begin{equation}\label{directsum1}
    {J_1}^{(r_1-r_2)}  \oplus {J_2}^{(r_2-r_3)} \oplus  \cdots \oplus J_{t-1}^{(r_{t-1}-r_{t})} \oplus
    J_t^{(r_{t})} \oplus A_t
 \end{equation}
 in which $J_k:=J_k(0)$
and $A_t$ is determined by $A$  up to
consimilarity and the other summands
are uniquely determined.
\end{theorem}
\proof
Let ${\cal A}: V
\to V$ be  a semilinear operator whose
matrix in some basis is $A$. Let
$W:={\cal A}V$ be the image of ${\cal
A}$. Then the matrix of the restriction
${\cal A}_1:W\to W$ of ${\cal A}$ on
$W$ is $A_1$. Applying Algorithm 1 to
$A_1$, we get the sequence $r_2, \dots,
r_{t}, A_t$. Reasoning by induction on
the length $t$ of the algorithm, we
suppose that  $r_2\ge r_3\ge \dots\ge
r_{t}$ and that $A_1$ is consimilar to
 \begin{equation}\label{directsum2}
    {J_1}^{(r_2-r_3)} \oplus  \cdots \oplus J_{t-2}^{(r_{t-1}-r_t)}\oplus
    J_{t-1}^{(r_{t})} \oplus A_t.
 \end{equation}
Thus, ${\cal A}_1:W\to W$ is given by
the matrix \eqref{directsum2} in some
basis of $W$.

The direct sum \eqref{directsum2}
defines the decomposition of $W$ into
the direct sum of invariant subspaces
\begin{equation*}
W=(W_{21}\oplus\dots\oplus
W_{2,r_2-r_3})\oplus\cdots \oplus
(W_{t1}\oplus\dots\oplus
W_{tr_t})\oplus W'.
\end{equation*}
Each $W_{pq}$ is generated by some
basis vectors $e_{pq2},$
$e_{pq3},\dots,e_{pqp}$ such that
\begin{equation*}\label{gr}
\mathcal A:\ e_{pq2}\mapsto e_{pq3}\mapsto \cdots\mapsto e_{pqp}\mapsto 0.
\end{equation*}

For each $W_{pq}$, we choose
$e_{pq1}\in V$ such that $\mathcal
Ae_{pq1}=e_{pq2}$. The set
$$\{e_{pqp}\,|\,2\le p\le t,\ 1\le q\le r_p-r_{p+1}\}\quad (r_{t+1}:=0)$$ consists of $r_2$ basis vectors belonging to the kernel of $\cal A$; we supplement this set to a basis of the kernel of $\cal A$ by some vectors $e_{111},\dots,e_{1,r_1-r_2,1}$.

The set of vectors $e_{pqs}$
supplemented by the vectors of some
basis of $W'$ is a basis of $V$. The
matrix of $\cal A$ in this basis has
the form \eqref{directsum1} because
\begin{equation*}\label{gri}
\mathcal A:\ e_{pq1}\mapsto e_{pq2}\mapsto e_{pq3}\mapsto \cdots\mapsto e_{pqp}\mapsto 0
\end{equation*}
for all $p=1,\dots,t$ and
$q=1,\dots,r_p-r_{p+1}$. This completes
the proof of Theorem \ref{thrm1}.

\begin{example}
Let a square
matrix $A$ define a semilinear operator
$\mathcal A:V\to V$ and let the
singular part of its regularizing
decomposition be ${J_2} \oplus J_{3}
\oplus J_4$. This means that $V$
possesses a set of linear independent
vectors forming the Jordan chains
\begin{align}\nonumber
\mathcal A:\quad&e_1\mapsto e_2\mapsto e_3\mapsto e_4\mapsto 0
\\\label{der}&f_1\mapsto f_2\mapsto f_3\mapsto 0
\\ \nonumber &g_1\mapsto g_2\mapsto 0
\end{align}
Applying the first step of Algorithm 1,
we get $A_1$ whose singular part
corresponds to the chains
\begin{align*}
\mathcal A:\quad&e_2\mapsto e_3\mapsto e_4\mapsto 0
\\&f_2\mapsto f_3\mapsto 0
\\&g_2\mapsto 0
\end{align*}
On the second step, we delete
$e_2,f_2,g_2$ and so on. Thus, $r_i$ is
the number of vectors in the $i$th
column of \eqref{der}: $r_1=3,\ r_2=3,\
r_3=2, r_4 =1$. We get the singular
part of regularizing decomposition of
$A$:
 \begin{equation*}
   {J_1}^{(r_1-r_2)}  \oplus  \cdots \oplus J_{t-1}^{(r_{t-1}-r_{t})} \oplus
    J_t^{(r_{t})} =
    {J_1}^{(3-3)}  \oplus {J_2}^{(3-2)} \oplus J_{3}^{(2-1)} \oplus
    J_4^{(1)} = {J_2} \oplus J_{3} \oplus
    J_4.
 \end{equation*}
 In particular, if
\begin{equation}\label{dahr}
A=
\begin{MAT}(b){cccccccccc}
0&0 &0 &0&&&&&&\scriptstyle \it e_1\\
1&0 &0 &0&&&&&&\scriptstyle \it e_2\\
0&1 &0 &0&&&&&&\scriptstyle \it e_3\\
0 &0&1& 0&&&&&&\scriptstyle \it e_4\\
 &&& &0 & 0&0&&&\scriptstyle \it f_1\\
 &&& &1 & 0&0&&&\scriptstyle \it f_2\\
 &&& &0 & 1&0&&&\scriptstyle \it f_3\\
 &&& & & &&0&0&\scriptstyle \it g_1\\
 &&& & & &&1&0&\scriptstyle \it g_2\\
\scriptstyle \it e_1&\scriptstyle \it e_2
&\scriptstyle \it e_3&\scriptstyle \it e_4
&\scriptstyle \it f_1 &\scriptstyle \it f_2
&\scriptstyle \it f_3&\scriptstyle \it g_1 &\scriptstyle \it g_2&
\addpath{(0,1,4)uuuuuuuuurrrrrrrrrdddddddddlllllllll}
\addpath{(4,1,0)uuuuuuuuu}
\addpath{(7,1,0)uuuuuuuuu}
\addpath{(0,3,0)rrrrrrrrr}
\addpath{(0,6,0)rrrrrrrrr}
\\
\end{MAT}\,,
\end{equation}
then we can apply Algorithm 1 using
only transformations of permutational
similarity and obtain
\begin{equation*}\label{daJ}
\begin{MAT}(b){cccccccccc}
0& 0& 0&&&&&&&\scriptstyle \it e_1\\
0& 0&0&&&&&&&\scriptstyle \it f_1\\
0& 0&0&&&&&&&\scriptstyle \it g_1\\
1&0&0&0&0&0&&&&\scriptstyle \it e_2\\
0&1&0&0&0&0&&&&\scriptstyle \it f_2\\
0&0&1&0&0&0&&&&\scriptstyle \it g_2\\
& & &1&0&0&0&0&&\scriptstyle \it e_3\\
 & & &0&1&0&0&0&&\scriptstyle \it f_3\\
 & & & & &&1&0&0&\scriptstyle \it e_4\\
\scriptstyle \it e_1&\scriptstyle \it f_1
&\scriptstyle \it g_1&\scriptstyle \it e_2
&\scriptstyle \it f_2 &\scriptstyle \it g_2
&\scriptstyle \it e_3&\scriptstyle \it f_3 &\scriptstyle \it e_4&
\addpath{(0,1,4)uuuuuuuuurrrrrrrrrdddddddddlllllllll}
\addpath{(3,1,0)uuuuuuuuu}
\addpath{(6,1,0)uuuuuuuuu}
\addpath{(8,1,0)uuuuuuuuu}
\addpath{(0,2,0)rrrrrrrrr}
\addpath{(0,4,0)rrrrrrrrr}
\addpath{(0,7,0)rrrrrrrrr}
\addpath{(3,1,4)uuuuuurrrrrr}
\addpath{(6,1,4)uuurrr}
\addpath{(8,1,4)ur}
\\
\end{MAT}\,
\end{equation*}
(all unspecified blocks are zero),
which is the Weyr canonical form of
\eqref{dahr}, see \cite{j-h}.
\end{example}

\subsection{Regularizing  unitary algorithm for matrix pairs under mixed equivalence}

\quad\\

We say that pairs of $m \times n$
matrices $(A,B)$ and $(A',B')$  are
\emph{mixed equivalent} if there exist
nonsingular $S$ and $R$ such that
\begin{equation*}\label{consim2}
     (SA{R},SB\bar{R})=(A',B').
\end{equation*}
The {\it direct sum} of
      matrix pairs $(A,B)$ and $(C,D)$ is defined as follows:
       \[
      {(A,B)\oplus(C,D)}=\left(\begin{bmatrix}
  A & 0 \\
  0 & C
\end{bmatrix},\: \begin{bmatrix}
  B & 0 \\
  0 & D
\end{bmatrix}\right).
\]
The canonical form of a matrix pair
under mixed equivalence was obtained by
Djokovi\'{c} \cite{djok} (his result
was extended to undirected cycles of
linear and semilinear mappings in
\cite{deb2}): \\
{\it Each pair $(A,B)$ of
matrices of the same size is mixed
equivalent to a direct sum, determined
uniquely up to permutation of summands,
of pairs of the following types:
\[
(I_n,J_n(\lambda )),\ (I_n,
\left(\begin{smallmatrix} 0&1\\\mu
&0\end{smallmatrix}\right)),\
 (J_n(0),I_n),\
(F_n,G_n),\ (F_n^T,G_n^T),
\]
in which ${\lambda\ge 0}$ and
${\mu\notin\mathbb R}$ or ${\mu<0}$.}

Thus, $(A,B)$ is mixed equivalent to a
direct sum of a pair $(\underline A,
\underline B)$ of nonsingular matrices
and summands of the types:
\[
(I_n,J_n(0)),\ (J_n(0),I_n),\
(F_n,G_n),\ (F_n^T,G_n^T),
\]
in which $(\underline A,\underline B)$
is determined up to mixed equivalence
and the other summands are uniquely
determined up to permutation. This sum
is called a \emph{regularizing
decomposition} of $(A,B)$. The
following algorithm admits to construct
a regularizing decomposition using only
unitary transformations.

\begin{algorithm}
{\it Let $(A,B)$ be a pair of matrices
of the same size in which the rows of
$A$ are linearly dependent. By unitary
transformations of rows, we reduce $A$
to the form
\[ S_1A=
\begin{bmatrix}
  0 \\
  A' \\
\end{bmatrix},
\quad S_1 \text{ is unitary, }
\]
in which the rows of $A'$ are linearly
independent. These transformations
change $B$:
\begin{equation*}\label{B1}
S_1B=
\begin{bmatrix}
 B' \\
 B'' \\
\end{bmatrix}.
\end{equation*}

By unitary transformations of columns,
we reduce $B'$ to the form $[B_{1}' \
0]$ in which the columns of $B_{1}'$
are linearly independent, and obtain
\begin{equation*}\label{B2}
B{R_1}=
\begin{bmatrix}
  B_{1}' & 0 \\
  \star & B_{1} \\
\end{bmatrix}, \quad R_1 \text{ is unitary}.
\end{equation*}
These transformations change $A$:
\begin{equation*}\label{AB1}
S_1A\bar{R_1}=
\begin{bmatrix}
 0_{k_1l_1} & 0 \\
  \star & A_{1} \\
\end{bmatrix}.
\end{equation*}
We apply the same procedure to
$(A_{1},B_1)$ and obtain
\begin{equation*}\label{A2B2}
(S_2A_1\bar{R_2},S_2B_1{R_2})=
\left(
\begin{bmatrix}
  0_{k_2l_2} & 0 \\
  \star & A_2 \\
\end{bmatrix},
\begin{bmatrix}
  B_{2}' & 0 \\
  \star & B_2 \\
\end{bmatrix}
\right),
\end{equation*}
in which the rows of $[\star \ A_2]$
are linearly independent, $S_2$ and
$R_2$ are unitary, and the columns of
$B_{2}'$ are linearly independent.

We repeat this procedure until we
obtain
\begin{equation*}\label{AtBt}
({S}_{t}A_{t-1}{\bar{R}}_{t},{S}_{t}B_{t-1}{{R}}_{t})
=
\left(
\begin{bmatrix}
  0_{k_tl_t} & 0 \\
  \star & A_t \\
\end{bmatrix},
\begin{bmatrix}
  B_{t}' & 0 \\
  \star & B_t \\
\end{bmatrix}
\right),
\end{equation*}
in which  the rows of $A_t$ are are
linearly independent. The result of the
algorithm is the sequence
\begin{equation*}\label{asx}
(k_1,l_1),\ (k_2,l_2),\ \dots,\ (k_t,l_t),\
(A_t,B_t).
\end{equation*}
}
\end{algorithm}
 For a matrix pair $(A,B)$ and a
nonnegative integer $n$, we write
\[
(A,B)^{(n)}:=
           \begin{cases}
             (0_{00},0_{00}), & \hbox{if $n=0$;} \\
\underbrace{(A,B)\oplus\dots\oplus (A,B)}
_{n\text{ summands}},
 & \hbox{if $n\ge 1$.}
           \end{cases}
\]
\begin{theorem}\label{thrm3}
Let $(A,B)$ be a pair of complex
matrices of the same size. Let us apply
Algorithm 2 to $(A,B)$ and obtain
$$(k_1,l_1),\ (k_2,l_2),\ \dots,\
(k_t,l_t),\
(A_t,B_t).$$ Let us apply Algorithm 2
to $(\p A,\p B):=(B_t^T,A_t^T)$ and
obtain
$$(\p k_1,\p l_1),\ (\p k_2,\p l_2),
\ \dots,\ (\p k_{\p t},\p l_{\p t}),\ (\p
A_{\p t},\p B_{\p t}).$$ Then $(A,B)$ is
mixed equivalent to
\begin{align*}
&(F_1,G_1)^{(k_1-l_1)}\oplus\dots\oplus
(F_{t-1},G_{t-1})^{(k_{t-1}-l_{t-1})}
 \oplus (F_t,G_t)^{(k_t-l_{t})}
        \\\oplus
&{(J_1,I_1)}^{(l_1-k_2)}  \oplus\dots\oplus
{(J_{t-1},I_{t-1})}^{(l_{t-1}-k_{t})}\oplus
{(J_t,I_t)}^{(l_{t})}
        \\\oplus
&(F_1^T,G_1^T)^{(\p k_1-\p l_1)}\oplus\dots\oplus
(F_{\p t-1}^T,G_{\p t-1}^T)^{(\p k_{\p t-1}-l_{\p t-1})}
 \oplus  (F_{\p t}^T,G_{\p t}^T)^{(k_{\p t}-l_{\p t})}
        \\\oplus
&{(I_1,J_1)}^{(\p l_1-\p k_2)}  \oplus\dots\oplus
{(I_{\p t-1},J_{\p t-1})}^{(\p l_{\p t-1}-\p k_{\p t})}\oplus
{(I_{\p t},J_{\p t})}^{(\p l_{\p t})}
        \\\oplus
&(\p B^T_{\p t},\p A^T_{\p t})
\end{align*}
(all exponents in parentheses are
nonnegative). The pair $(\p B^T_{\p
t},\p A^T_{\p t})$ consists of
nonsingular matrices; it is determined
up to mixed equivalence. The other
summands are uniquely determined by
$(A,B)$.
\end{theorem}
The rows of $A_t$ in Theorem
\ref{thrm3} are linearly independent,
and so the columns of $\p B:=A_t^T$ are
linearly independent. As follows from
Algorithm 2, the columns of $\p B_{\p
t}$ are linearly independent too. Since
the rows of $\p A_{\p t}$ are linearly
independent, the columns of $\p B_{\p
t}$ are linearly independent, and the
matrices in $(\p A_{\p t},\p B_{\p t})$
have the same size, these matrices are
square, and so they are nonsingular.
The pairs $(I_n,J_n^T)$ and
$(G_n^T,F_n^T)$ are permutationally
equivalent to $(I_n,J_n)$ and
$(F_n^T,G_n^T)$.  Therefore, Theorem
\ref{thrm3} follows from the following
lemma.

\begin{lemma}\label{thrm2}
Let $(A,B)$ be a pair of complex
matrices of the same size. Let us apply
Algorithm 2 to $(A,B)$ and obtain
$$(k_1,l_1),\ (k_2,l_2),\ \dots,\
(k_t,l_t),\
(A_t,B_t).$$ Then $(A,B)$ is mixed
equivalent to
\begin{align}\nonumber
&(F_1,G_1)^{(k_1-l_1)}\oplus\dots\oplus
(F_{t-1},G_{t-1})^{(k_{t-1}-l_{t-1})}
 \oplus (F_{t},G_{t})^{(k_t-l_{t})}
        \\ \label{bgf}
        \oplus
&{(J_1,I_1)}^{(l_1-k_2)}  \oplus\dots\oplus
{(J_{t-1},I_{t-1})}^{(l_{t-1}-k_{t})}
        \\ \nonumber
    \oplus
&{(J_t,I_t)}^{(l_{t})}    \oplus
(A_{t},B_{t})
\end{align}
(all exponents in parentheses are
nonnegative). The rows of $A_t$ are
linearly independent. The pair
$(A_{t},B_{t})$  is determined up to
mixed equivalence. The other summands
are uniquely determined by $(A,B)$.
\end{lemma}

\proof
We write
\[
(A,B)\Longrightarrow (k_1,l_1,(A_1,B_1))
\]
if $k_1,l_1,(A_1,B_1)$ are obtained
from $(A,B)$ in the first step of
Algorithm 2.

First we prove two statements.

{\it Statement 1: If
\begin{equation}\label{eee}
\begin{aligned}
(A,B)\Longrightarrow
(k_1,l_1,(A_1,B_1)),
      \\
(\widetilde{A},\widetilde{B})\Longrightarrow
(\tilde{k}_1,\tilde{l}_1,(\widetilde{A}_1,\widetilde{B}_1)),
\end{aligned}
\end{equation}
and $(A,B)$ is mixed equivalent to
$(\widetilde{A},\widetilde{B})$, then
$k_1=\tilde{k}_1,$ $l_1= \tilde{l}_1$,
and $(A_1,B_1)$ is mixed equivalent to
$(\widetilde{A}_1,\widetilde{B}_1)$. }

Let $m$ be the number of rows in $A$.
Then
\[
k_1=m-\rank A=m-\rank \widetilde{A}=\tilde{k}_1.
\]
Since $(A,B)$ and
$(\widetilde{A},\widetilde{B})$ are
mixed equivalent and they are reduced
by mixed equivalence transformations to
\begin{equation}\label{mkj}
\setlength{\arraycolsep}{2pt}
\left(
    \begin{bmatrix}
  0_{k_1l_1} & 0 \\
  X & A_{1} \\
 \end{bmatrix},
  \begin{bmatrix}
  B_1' & 0 \\
  Y & B_{1} \\
 \end{bmatrix}
 \right),
 \left(
 \begin{bmatrix}
  0_{{k}_1\tilde{l}_1} & 0 \\
  \widetilde{X} & \widetilde{A}_1 \\
 \end{bmatrix},
  \begin{bmatrix}
  \widetilde{B}_1' & 0 \\
  \widetilde{Y} & \widetilde{B}_1 \\
 \end{bmatrix}
 \right),
\end{equation}
there exist nonsingular $S$ and $R$
such that
\begin{equation}\label{pfprop}
\renewcommand\arraystretch{1.2}
\left(S
    \begin{bmatrix}
  0_{k_1l_1} & 0 \\
  X & A_{1} \\
 \end{bmatrix},\ S
  \begin{bmatrix}
  B_1' & 0 \\
  Y & B_{1} \\
 \end{bmatrix}
 \right)\\
 =\left(
 \begin{bmatrix}
  0_{{k}_1\tilde{l}_1} & 0 \\
  \widetilde{X} & \widetilde{A}_1 \\
 \end{bmatrix} R,\
  \begin{bmatrix}
  \widetilde{B}_1' & 0 \\
  \widetilde{Y} & \widetilde{B}_1 \\
 \end{bmatrix}\bar{R}
 \right).
\end{equation}
Equating the first matrices of these
pairs, we find that $S$ has the form
 \begin{equation*}\label{eq5}
S= \begin{bmatrix}
  S_{11}&0  \\
  S_{21}&S_{22} \\
 \end{bmatrix},\quad S_{11}\text{ is }k_1\times k_1.
  \end{equation*}
Equating the second matrices of the
pairs \eqref{pfprop}, we find that
\begin{equation}\label{nnm}
  S_{11}[B_1'\ 0]=
  [\widetilde{B}_1' \ 0]
\bar{R},
\end{equation}
and so \[l_1=\rank [B_1' \ 0]=\rank
[\widetilde{B}_1' \ 0]=\tilde{l}_1.\]
Since $B_1'$ and $\widetilde{B}_1'$ are
$k_1\times l_1$ and have linearly
independent columns, \eqref{nnm}
implies that $R$ is of the form
\begin{equation*}
R=\begin{bmatrix}
R_{11}&0\\
R_{21}&R_{22}
\end{bmatrix},\quad R_{11}\text{ is }l_1\times l_1.
\end{equation*}
Equating the (2,2) entries in the
matrices \eqref{pfprop}, we get
\[
S_{22}A_1= \widetilde{A}_1R_{22},\qquad
S_{22}B_1= \widetilde{B}_1\bar{R}_{22},
\]
hence $(A_1,B_1)$ and
$(\widetilde{A}_1,\widetilde{B}_1)$ are
mixed equivalent, which completes the
proof of Statement 1.
\medskip

{\it Statement 2: If \eqref{eee}, then
\[
(A,B)\oplus (\widetilde{A},\widetilde{B})  \Longrightarrow
     (k_1+\tilde{k}_1,  l_1+\tilde{l}_1,(A_1 \oplus \widetilde{A}_1 ,B_1 \oplus \widetilde{B}_1)).
\]
} Indeed, if $(A,B)$ and
$(\widetilde{A},\widetilde{B})$ are
reduced to \eqref{mkj}, then
$(A,B)\oplus
(\widetilde{A},\widetilde{B})$ is
reduced to
  \begin{equation*}
 \renewcommand\arraystretch{1.2}
 \left(
     \begin{bmatrix}
 0_{k_1l_1}\oplus 0_{\tilde{k}_1\tilde{l}_1}
   & 0 \oplus 0 \\
 X\oplus \widetilde{X} &
 A_1\oplus\widetilde{A}_1
 \end{bmatrix},
  \begin{bmatrix}
 B_1'\oplus \widetilde{B}_1' & 0\oplus 0  \\
  Y\oplus \widetilde{Y} &
  B_1\oplus \widetilde{B}_1
 \end{bmatrix}
 \right),
 \end{equation*}
which is permutationally equivalent to
\[
\renewcommand\arraystretch{1.2}
\setlength{\arraycolsep}{2pt}
\left(
    \begin{bmatrix}
  0_{k_1l_1} & 0 \\
  X & A_{1} \\
 \end{bmatrix}\oplus
  \begin{bmatrix}
  B_1' & 0 \\
  Y & B_{1} \\
 \end{bmatrix}
 \right),\
\left(
 \begin{bmatrix}
  0_{\tilde{k}_1\tilde{l}_1} & 0 \\
   \widetilde{X}& \widetilde{A}_1 \\
 \end{bmatrix}\oplus
  \begin{bmatrix}
  \widetilde{B_1'} & 0 \\
 \widetilde{Y} & \widetilde{B}_1 \\
 \end{bmatrix}
 \right).
\]

We are ready to prove Lemma \ref{thrm2}
for any pair $(A,B)$. Due to Statement
1, we can replace $(A,B)$ by any mixed
equivalent pair. In particular, we can
take
\begin{align}\label{tfe}
(A,B)=&(F_{1},G_{1})^{(r_1)}\oplus\dots\oplus
 \ca (F_{t},G_{t})^{(r_t)}
        \oplus \\ \nonumber
&{(J_1,I_1)}^{(s_1)}  \oplus\dots\oplus
{(J_t,I_t)}^{(s_{t})}
\oplus
(C,D)
\end{align}
for some nonnegative
$t,r_1,\dots,r_t,s_1,\dots,r_t$ and
some pair $(C,D)$ in which $C$ has
linearly independent rows.

Clearly,
\[
    (J_i,I_i)\Longrightarrow
\begin{cases}
(1,1,(J_{i-1},I_{i-1})), &\text{if }i \neq 1; \\
(1,1,(0_{00},0_{00})), &\text{if } i=1,
\end{cases}
\]
and
\[
    (F_{i},G_{i})\Longrightarrow
    \begin{cases}
(1,1,(F_{i-1},G_{i-1})), &\text{if } i \neq 1; \\
(1,0,(0_{00},0_{00})), &\text{if } i=1.
\end{cases}
\]

Due to Statement 2,
\begin{itemize}
  \item $k_1=m-\rank A$ is the
      number of all summands of the
      types $(J_i,I_i)$ and $(F_{i},G_{i})$,

  \item $l_1$ is the number of all
      summands of the types
      $(J_i,I_i)$ and $(F_{i},G_{i})$,
      except for $(F_{1},G_{1})$,

  \item and
\begin{align}\label{tfe1}
(A_1,B_1)=&(F_{1},G_{1})^{(r_2)}\oplus\dots\oplus
 (F_{t-1},G_{t-1})^{(r_t)}
        \oplus \\ \nonumber
&{(J_1,I_1)}^{(s_2)}  \oplus\dots\oplus
{(J_{t-1},I_{t-1})}^{(s_{t})}
        \oplus
(C,D).
\end{align}
\end{itemize}
We find that $k_1-l_1$ is the number of
summands of the type $(F_{1},G_{1})$.

Applying the same reasoning to
\eqref{tfe1} instead of \eqref{tfe} we
get that
\begin{itemize}
  \item $k_2$ is the number of all
      summands of the types
      $(J_i,I_i)$ and $(F_{i},G_{i})$
      with $i\ge 2$,

  \item $l_1$ is the number of all
      summands of the types
      $(J_i,I_i)$ with $i\ge 2$ and
      $(F_{i},G_{i})$ with $i\ge 3$,

  \item and
\begin{align*}\label{tfeu}
(A_2,B_2)=&(F_{1},G_{1})^{(r_3)}\oplus\dots\oplus
 (F_{t-2},G_{t-2})^{(r_t)}
        \oplus\\ \nonumber
&{(J_1,I_1)}^{(s_3)}  \oplus\dots\oplus
{(J_{t-2},I_{t-2})}^{(s_{t})}
        \oplus
(C,D).
\end{align*}
\end{itemize}
We find that $k_2-l_2$ is the number of
summands of the type $(F_{1},G_{1})$, and
that $l_1-k_2$ is the number of
summands of the type $(J_1,I_1)$, and
so on, until we obtain \eqref{bgf}.

The fact that the pair $(A_{t},B_{t})$
in \eqref{bgf} is determined up to
mixed equivalence and the other
summands are uniquely determined by
$(A,B)$ follows from Statement 1 (or
from the canonical form of a matrix
pair up to mixed equivalence). This
concludes the proof of Lemma
\ref{thrm2} and Theorem 2.


\begin{thebibliography}{99}

\bibitem{deb2}
    D.D. de Oliveira, R.A. Horn,
    T. Klimchuk, V.V. Sergeichuk, (2012),
    Remarks
    on the classification of a pair of
    commuting semilinear operators,
    Linear Algebra Appl. 436
    3362--3372.
    
\bibitem{djok} D.\v{Z}.
    Djokovi\'{c}, (1978),
    Classification of pairs consisting
    of a linear and a semilinear map,
    Linear Algebra Appl. 20
    147--165.

\bibitem{hongh} Y.P. Hong, R.A.
      Horn, (1988), A canonical form
      for matrices under
      consimilarity, Linear
      Algebra Appl. 102
      143--168.

\bibitem{hornj} R.A. Horn, C.R.
      Johnson, (2012), Matrix
      Analysis, 2nd ed., Cambridge
      University Press, New York.

\bibitem{h-s} R.A. Horn, V.V.
    Sergeichuk, (2006),
A regularization algorithm for matrices
of bilinear and sesquilinear forms,
Linear Algebra Appl. 412
380--395.



\bibitem{ser} V.V. Sergeichuk, (2004),
 Computation of
    canonical matrices for chains and
    cycles of linear mappings,
Linear Algebra Appl. 376
235--263.

\bibitem{van}
P. Van Dooren, (1979),
The computation of Kronecker's canonical form of a singular pencil, Linear Algebra
Appl. 27  103--140.

\bibitem{var} A. Varga, (2004),
Computation of Kronecker-like forms of periodic matrix pairs, Symp. on Mathematical Theory of Networks and Systems, Leuven, Belgium, July 5--9.






\end{thebibliography}
\end{document}